\documentclass[12pt,letterpaper,reqno]{amsart}
\usepackage[T1]{fontenc}
\usepackage{lmodern}
\usepackage[hmargin=1in,vmargin=0.9in]{geometry}
\usepackage{amsmath,amssymb,mathtools}
\usepackage[expansion=false]{microtype}
\usepackage{xcolor,graphicx,booktabs}
\usepackage{flafter}
\usepackage[colorlinks=true,linkcolor=blue!75!black,citecolor=red!75!black,urlcolor=blue!60!black]{hyperref}

\newtheorem{theorem}{Theorem}[section]
\newtheorem{lemma}[theorem]{Lemma}

\newtheorem{corollary}[theorem]{Corollary}
\numberwithin{equation}{section}
\newcommand{\R}{\mathbb R}
\newcommand{\C}{\mathbb C}
\newcommand{\E}{\mathbb E}

\newcommand{\HS}{\mathrm{HS}}
\newcommand{\op}{\mathrm{op}}
\newcommand{\sinc}{\operatorname{sinc}}
\newcommand{\norm}[1]{\left\lVert#1\right\rVert}
\newcommand{\abs}[1]{\left\lvert#1\right\rvert}
\newcommand{\ip}[2]{\left\langle#1,#2\right\rangle}
\allowdisplaybreaks[2]
\title[Optimal approximation of stochastic fractional Helmholtz fields]{Optimal approximation of stochastic fractional Helmholtz fields}
\hypersetup{pdftitle={Optimal approximation of stochastic fractional Helmholtz fields}}
\author[D. Qiu]{Dong Qiu}
  \address{School of Mathematical Sciences, Zhejiang University, Hangzhou 310058, China}
\email{qiudong@zju.edu.cn}

  \author[X. Xu]{Xiang Xu}
  \address{School of Mathematical Sciences, and Center for Interdisciplinary Applied Mathematics, Zhejiang University, Hangzhou 310058, China}
\email{xxu@zju.edu.cn}

 \author[Y. Zhao]{Yue Zhao}
\address{School of Mathematics and Statistics, and Key Lab NAA-MOE, Central China Normal University,
Wuhan 430079, China}
\email{zhaoy@ccnu.edu.cn}

\date{}

\begin{document}
\begin{abstract}
We consider the outgoing solution of the fractional Helmholtz equation driven by spatial white noise with prescribed source strength in two and three dimensions. The data are the integrals of a single realization of the white noise over the cells of nested dyadic meshes. We divide the squared error integrated over the observation region by a power of the source mesh size. The exponent is four times the fractional power of the Laplacian minus the spatial dimension. For each fixed sequence of approximations based on the cell integrals, the limit inferior of this rescaled error is bounded below by an explicit constant almost surely. The piecewise constant approximation of the white noise gives the conditional mean, which attains this constant and has the same leading constant in mean square. The measures defined by its rescaled squared error converge weakly almost surely to a measure with density proportional to the source strength on the source region. On aligned boxes, we construct a cell-averaged fast Fourier transform method that attains the optimal constant when the ratio of the output mesh size to the source mesh size tends to zero. We test the predicted constants by comparing successive approximations on the same realization.
\end{abstract}
\maketitle

\section{Introduction}

Stochastic partial differential equations describe physical systems subject to uncertain forcing. In wave propagation, spatial white noise is used to model fluctuations in the source. We study the accuracy with which the radiated field can be computed for a given realization of the noise.

Numerical methods for stochastic parabolic equations include lattice and time discretizations~\cite{Gyongy,JK,MGR}. The convergence of finite difference, trigonometric, and Verlet schemes has been analyzed for stochastic wave equations~\cite{Walsh,CLS,BLM}. For the fractional Helmholtz equation, numerical methods have been developed for the deterministic problem~\cite{Glusa}, and inverse source problems have been studied in homogeneous and inhomogeneous media~\cite{LiLi,QXZinhom,QXZfree}. Related work on fractional elliptic equations includes mean square error estimates for equations driven by spatial white noise and the approximation of spherical Gaussian random fields~\cite{BKK,JKL}.

We consider the outgoing solution of the stochastic fractional Helmholtz equation
\begin{equation*}
 \bigl((-\Delta)^\alpha-k^{2\alpha}\bigr)u
   =\sqrt{\mu(x)}\,\dot W_D
   \qquad\text{in }\R^d,\quad d=2,3 .
\end{equation*}
Here $k>0$ is a fixed wave number, $d/4<\alpha<1$, and $\dot W_D$ is spatial white noise supported in a bounded source region $D$. The source strength function $\mu$ is prescribed, and $\sqrt\mu$ is the coefficient multiplying the white noise. We approximate the solution $u$ in $L^2(Q)$ on a bounded observation region $Q$ containing $\overline D$. The data are the integrals of the underlying white noise over Cartesian mesh cells. The meshes are nested, and the cell integrals at every mesh size come from the same realization.

\begin{samepage}
For the classical stochastic Helmholtz equation, Cao, Zhang, and Zhang~\cite[Theorem~1]{CZZ} derived a mean square error estimate for the approximate solution obtained by replacing the white noise with its piecewise constant approximation. Finite difference methods have also been studied for the one-dimensional stochastic Helmholtz equation~\cite{CTZ}. For the fractional Helmholtz equation, convergence of piecewise constant approximations of the noise has been used to construct the outgoing solution~\cite{QXZfree}. We seek the best accuracy attainable from the cell integrals on a single realization.\par
\end{samepage}

We establish the almost-sure optimality of the conditional mean for approximating the outgoing solution from the cell integrals of a single realization. Let $h$ denote the source mesh size. For each fixed sequence of approximations based on the cell integrals, the limit inferior of the squared $L^2(Q)$ error divided by $h^{4\alpha-d}$ is at least an explicit constant $C_\mu$ almost surely along the dyadic meshes. For the conditional mean, the rescaled squared error converges almost surely to $C_\mu$, the same constant as in the mean square limit. The constant is proportional to the total source strength $\int_D\mu(x)\,dx$. The measures defined by the rescaled squared error of the conditional mean converge weakly almost surely to a measure with density proportional to $\mu(x)$ on $D$.

On aligned boxes, we construct a cell-averaged fast Fourier transform (FFT) method and determine how the output mesh affects the leading error constant. The method attains $C_\mu$ when the ratio of the output mesh size to the source mesh size tends to zero. We derive the limits for differences between successive approximations on separate output grids and on a common finer grid. We test these limits with a single realization without evaluating the exact solution.

Section~\ref{sec:model} introduces the model and the approximate problem. Sections~\ref{sec:sharp} and~\ref{sec:solver} establish the optimal error estimates and analyze the FFT method, respectively. Section~\ref{sec:experiment} presents the numerical examples.

\section{The model and the approximate problem}\label{sec:model}

We work in two and three dimensions. Let $D$ and $Q$ be bounded open subsets of $\R^d$ with $D$ nonempty and $\overline D\subset Q$. We fix $k>0$ and $d/4<\alpha<1$, and set $p=4\alpha-d$. Thus $0<p<2$. On a probability space $(\Omega,\mathcal F,\mathbb P)$, let $W$ be real isonormal Gaussian white noise on $L^2(D)$, with $\E(W(v)W(w))=\int_D vw$ for real $v,w\in L^2(D)$. Let $\mu$ be the prescribed source strength function, and let $a=\sqrt\mu$ be the coefficient multiplying the white noise. We assume that $a$ is bounded and Lipschitz on $D$ and that the total source strength is $I_\mu=\int_D\mu(x)\,dx>0$. We define the random source by $\langle \sqrt\mu\,\dot W_D,\varphi\rangle=W(a\varphi|_D)$ for real test functions $\varphi$. Its covariance on test functions $\varphi,\psi$ is $\int_D\mu\varphi\psi$. We consider
\begin{equation}\label{eq:model}
  \bigl((-\Delta)^\alpha-k^{2\alpha}\bigr)u
   =\sqrt{\mu(x)}\,\dot W_D
   \quad\text{in }\R^d .
\end{equation}
We use the Fourier convention $\widehat v(\xi)=\int_{\R^d}e^{-ix\cdot\xi}v(x)\,dx$. The fractional Laplacian is defined by
\begin{equation*}
  \widehat{(-\Delta)^\alpha v}(\xi)=\abs{\xi}^{2\alpha}\widehat v(\xi).
\end{equation*}
The outgoing Green function $G$ is defined by the limiting absorption prescription~\cite{ZCV}
\begin{equation}\label{eq:G}
  \widehat G(\xi)=\bigl(\abs{\xi}^{2\alpha}-k^{2\alpha}-i0\bigr)^{-1}.
\end{equation}
Here $-i0$ means that the denominator is shifted by $-i\varepsilon$ with $\varepsilon>0$ and the limit $\varepsilon\downarrow0$ is taken in the sense of tempered distributions. The negative imaginary shift selects the Green function satisfying the outgoing radiation condition at infinity.

Let $T$ be the integral operator with kernel $G(x-y)$ on $Q\times D$. Write $M_a$ for multiplication by $a$ on $D$ and put $T_\mu=TM_a$. We denote the Hilbert--Schmidt and operator norms by $\norm{\cdot}_{\HS}$ and $\norm{\cdot}_{\op}$, respectively. For the Gaussian Hilbert space framework, see~\cite[Chapters~1 and~9]{Janson}. Given a Hilbert--Schmidt operator $A$ from the real input space into $L^2(Q;\C)$, write $AW$ for the mean square limit of $\sum_j W(e_j)Ae_j$ in a real orthonormal basis $(e_j)$. The white noise isometry $\E\norm{AW}^2=\norm A_{\HS}^2$ follows from orthogonality. The complex output space is viewed as a real Hilbert space, with inner product given by the real part of its complex inner product. We write $\norm{\cdot}$ for the $L^2(Q)$ norm unless another norm is specified. Constants in estimates may depend on the fixed model and regions and are independent of the mesh sizes. We use the same notation for $T$ with arbitrary bounded input and output sets in the following lemma.

\begin{lemma}\label{lem:riesz}
Let $K$ be the Riesz kernel with $\widehat K(\xi)=\abs{\xi}^{-2\alpha}$. Then $G=K+G_{\rm rem}$, where $G_{\rm rem}\in H^1_{\rm loc}(\R^d)$. The kernel $G(x-y)$ is square integrable on bounded products. The operators $T$ and $T^*$ map $L^2$ boundedly into $H^1$ on bounded sets, including bounded neighborhoods of their respective output sets.
\end{lemma}

\begin{proof}
We have $0<2\alpha<d$ and $K(x)=c_{d,\alpha}\abs x^{2\alpha-d}$ with $c_{d,\alpha}>0$ by the Fourier transform formula for the Riesz kernel~\cite[Chapter~V, Section~1]{Stein}. As in~\cite[Lemma~2.1]{QXZfree}, we split the Fourier transforms by a smooth cutoff equal to one on a neighborhood of the closed ball of radius $2k+1$. The inverse transforms of the low-frequency parts are smooth because their Fourier transforms have compact support. At high frequencies,
\begin{equation*}
 \bigl(\abs{\xi}^{2\alpha}-k^{2\alpha}\bigr)^{-1}-\abs{\xi}^{-2\alpha}
 =O\bigl(\abs{\xi}^{-4\alpha}\bigr).
\end{equation*}
Since $8\alpha>d+2$, the inverse transform of the high-frequency remainder belongs to $H^1(\R^d)$. The condition $4\alpha>d$ gives local square integrability of $K$. The high-frequency part of the multiplier in~\eqref{eq:G} maps $L^2$ into $H^{2\alpha}$, and $2\alpha>1$. The same argument applies to the adjoint, with conjugation and the real part on the real input space.
\end{proof}

In particular, $u=T_\mu W\in L^2(\Omega;L^2(Q))$. This is the outgoing distributional solution defined by~\eqref{eq:G}. For any $s>d/2$, the map taking $f\in L^2(D)$ to the zero extension of $af$ is Hilbert--Schmidt into $H^{-s}(\R^d)$. Its squared Hilbert--Schmidt norm is $(2\pi)^{-d}I_\mu\int_{\R^d}(1+|\xi|^2)^{-s}\,d\xi$ by Parseval's identity. This gives a compactly supported random distribution whose Fourier transform multiplies the outgoing boundary value. The factor $\abs{\xi}^{2\alpha}$ in the equation acts by ordinary multiplication near the origin, where the multiplier is locally bounded.

For a given source mesh size $h$, intersect the Cartesian cubes of side $h$ with $D$ and retain the cells of positive measure. Let $P_h$ denote the $L^2(D)$ orthogonal projection onto functions that are constant on each cell, and set $a_h=P_h a$. On each cell $E$, the piecewise constant approximation $P_hW$ of the white noise equals the cell average $W(\mathbf 1_E)/|E|$. For earlier piecewise constant approximations of white noise, see~\cite{ANZ} and~\cite[Section~2]{DZ}. This construction has been used for Helmholtz and semilinear elliptic problems~\cite{CZZ,CYI,QXZfree}. Cao, Hong, and Liu~\cite{CHL} studied a corresponding approximation for a boundary value problem driven by fractional Brownian motion.

The data at source mesh size $h$ are the cell integrals $W(\mathbf 1_E)$, and $\mathcal F_h$ denotes the sigma-algebra they generate. We write $u_h=T_\mu P_hW$ for the solution of the approximate problem obtained by replacing $\dot W_D$ by $P_hW$ in~\eqref{eq:model}. Here $P_hW$ is extended by zero outside $D$.

We use the nested dyadic meshes $h_m=h_0 2^{-m}$ with a common grid vertex. The white noise integral over each coarse cell equals the sum of the integrals over the fine cells it contains. This identity applies also to cells cut by the boundary of $D$. Thus $\mathcal F_{h_m}\subset\mathcal F_{h_{m+1}}$, and all the cell integrals come from the same realization of $W$.

\section{Optimal error estimates from the cell integrals}\label{sec:sharp}

The error $u-u_h$ is obtained by applying $T_\mu$ to $(I-P_h)W$, which is independent of the cell integrals. The singular part of the Green function in Lemma~\ref{lem:riesz} determines the leading term in the mean square $L^2$ error. The smoother remainder contributes a term of higher order.

Write $C_0=(-1/2,1/2)^d$ for the unit cube and $\sinc(t)=\sin(t)/t$, with its continuous value at zero. Define
\begin{equation}\label{eq:constant}
  \mathcal C_{d,\alpha}
  =\frac{1}{(2\pi)^d}\int_{\R^d}\abs{\xi}^{-4\alpha}
       \left(1-\prod_{\ell=1}^d\sinc^2(\xi_\ell/2)\right)\,d\xi,
  \qquad C_\mu=I_\mu\mathcal C_{d,\alpha}.
\end{equation}
The integral is finite and positive since the factor in parentheses vanishes quadratically at the origin and is bounded at infinity. For later numerical evaluation, an equivalent expression is
\begin{equation}\label{eq:constant-cell}
 \mathcal C_{d,\alpha}
 = b_{d,p}\int_{(-1,1)^d}\abs z^p\prod_{\ell=1}^d(1-\abs{z_\ell})\,dz,\qquad
 b_{d,p}=\frac{\Gamma(1-p/2)}{2^{d+p-1}\pi^{d/2}p\,\Gamma((d+p)/2)} .
\end{equation}
Indeed, the factor in parentheses in~\eqref{eq:constant} equals the integral of $1-\cos(\xi\cdot z)$ against $\prod_\ell(1-\abs{z_\ell})$ on $(-1,1)^d$. The identity
$\int_0^\infty t^{(d+p)/2-1}e^{-t\abs{\xi}^2}\,dt
=\Gamma((d+p)/2)\abs{\xi}^{-d-p}$ and the Gaussian Fourier integral give~\eqref{eq:constant-cell}.

The next theorem gives the weighted mean square error of $u_h$ for any nonnegative continuous function $\phi$. Setting $\phi=1$ gives the mean square $L^2(Q)$ error.

\begin{theorem}\label{thm:risk}
For every $h$, one has $u_h=\E(u\mid\mathcal F_h)$. For each square-integrable $\mathcal F_h$-measurable approximation $Y_h$,
\begin{equation}\label{eq:risk}
 \E\norm{u-Y_h}^2=R_h+\E\norm{u_h-Y_h}^2,
 \qquad R_h=\norm{T_\mu(I-P_h)}_{\HS}^2.
\end{equation}
For $\phi\in C(\overline Q)$ with $\phi\ge 0$,
\begin{equation}\label{eq:local-risk}
 h^{-p}\norm{\phi^{1/2}T_\mu(I-P_h)}_{\HS}^2
  \longrightarrow\mathcal C_{d,\alpha}\int_D\phi(x)\mu(x)\,dx .
\end{equation}
In particular, $R_h\sim C_\mu h^p$.
\end{theorem}

\begin{proof}
The projection and its orthogonal complement involve independent Gaussian coordinates. Hence $X_h=T_\mu(I-P_h)W$ is independent of $\mathcal F_h$ and centered. Conditional expectation and the isometry give~\eqref{eq:risk}.

Because $a$ is Lipschitz, $\norm{a-a_h}_{L^\infty(D)}\le Ch$, including on clipped cells. The difference between $\phi^{1/2}T_\mu(I-P_h)$ and $\phi^{1/2}T(I-P_h)M_{a_h}$ has Hilbert--Schmidt norm $O(h)$, since $T$ is Hilbert--Schmidt and $M_{a_h}$ commutes with $P_h$. As $p<2$, it suffices to find the leading error for this cellwise constant source coefficient.

On a clipped cell, the mean of a function minimizes its squared error over constants. If the function is defined on the full cube, its mean over that cube can be used as a comparison constant. Applying the cube Poincar\'e inequality to $G_{\rm rem}(x-\cdot)$ and integrating in $x$ gives
\begin{equation*}
 \norm{\phi^{1/2}T_{\rm rem}(I-P_h)M_{a_h}}_{\HS}\le C_\phi h ,
\end{equation*}
where $T_{\rm rem}$ has kernel $G_{\rm rem}(x-y)$ and the kernel gradient is integrated on a fixed bounded enlargement of $Q\times D$. Thus it suffices to compute the limit for $K$.

Set
\begin{equation*}
 g(t)=\int_{C_0}\left|K(t-z)-\int_{C_0}K(t-v)\,dv\right|^2 dz .
\end{equation*}
The cell variance identity and Plancherel's theorem give
\begin{equation*}
 \int_{\R^d}g(t)\,dt
 =\frac12\int_{C_0}\int_{C_0}\norm{K(\cdot-z)-K(\cdot-v)}_{L^2(\R^d)}^2\,dz\,dv
 =\mathcal C_{d,\alpha}.
\end{equation*}
Since $0<p<2$, the differences belong to $L^2(\R^d)$. By homogeneity,
\begin{equation*}
 \norm{K(\cdot-z)-K(\cdot-v)}_{L^2(\R^d)}^2\le C\abs{z-v}^p .
\end{equation*}
The variance identity and boundedness of $a_h$ therefore bound the contribution of the clipped cells by $C_\phi h^p$ times their total volume. Clipped cells lie in the interior strip $\{y\in D:\operatorname{dist}(y,\R^d\setminus D)\le\sqrt d\,h\}$, whose volume tends to zero.

For a full cube $E$ centered at $c$, let $a_E$ be the value of $a_h$ on $E$. Homogeneity gives its contribution as
\begin{equation*}
 a_E^2h^{p+d}\int_{\R^d}\widetilde\phi(c+ht)g(t)\,dt ,
\end{equation*}
where $\widetilde\phi$ is extended by zero outside $Q$. The difference between the integral and $\phi(c)\int g$ tends to zero uniformly for $c\in D$, by first using uniform continuity on bounded $t$ and then the integrable tail of $g$. Here $\overline D\subset Q$ guarantees that the bounded-$t$ portion lies in $Q$ for small $h$. The weighted Riemann sums over the full cubes converge to $\int_D\phi(x)\mu(x)\,dx$ because $a_h$ converges uniformly to $a$. The $O(h)$ bounds and $p<2$ complete the proof.
\end{proof}

The connection between average-case and almost-sure optimality for Gaussian linear problems was studied by~\cite{WW}. We use a variance estimate for the normalized squared error of $u_{h_m}$ to obtain its almost-sure limit along the dyadic meshes. Independence of $u-u_{h_m}$ from the cell integrals then gives an asymptotic error formula for each fixed sequence of approximations based on these data.

\begin{theorem}\label{thm:path}
Let $h_m=h_0 2^{-m}$. For every fixed sequence of $\mathcal F_{h_m}$-measurable $L^2(Q)$-valued approximations $Y_m$, one has almost surely
\begin{equation}\label{eq:path-pythagoras}
 \norm{u-Y_m}^2
 =C_\mu h_m^p+\norm{u_{h_m}-Y_m}^2+
 o\bigl(h_m^p+\norm{u_{h_m}-Y_m}^2\bigr).
\end{equation}
In particular, $\liminf_m h_m^{-p}\norm{u-Y_m}^2\ge C_\mu$, and the limit is $C_\mu$ precisely on the paths on which $\norm{u_{h_m}-Y_m}=o(h_m^{p/2})$. The exceptional null event may depend on the sequence. The same conclusions hold with additional randomness independent of $W$. Moreover,
\begin{equation}\label{eq:weak}
 h_m^{-p}\abs{u-u_{h_m}}^2\,dx
 \ \rightharpoonup\ \mathcal C_{d,\alpha}\mu\mathbf 1_D\,dx
 \quad\text{on }\overline Q\quad\text{almost surely}.
\end{equation}
\end{theorem}

\begin{proof}
Let $A_h=T_\mu(I-P_h)$. Applying the cube Poincar\'e inequality to the extended adjoint in Lemma~\ref{lem:riesz} gives $\norm{T(I-P_h)}_{\op}\le Ch$, using the full-cube comparison on clipped cells. The cellwise constant source coefficient commutes with $P_h$. Therefore
\begin{equation}\label{eq:op}
 \norm{A_h}_{\op}\le\norm{T(I-P_h)}_{\op}\norm{a_h}_{L^\infty(D)}+\norm T_{\op}\norm{a-a_h}_{L^\infty(D)}\le Ch .
\end{equation}
For any deterministic Hilbert--Schmidt operator $B$, expansion in the eigenvectors of its Gaussian covariance and the Gaussian fourth-moment identity~\cite{Isserlis} give
\begin{equation}\label{eq:variance}
 \operatorname{Var}\bigl(\norm{BW}^2\bigr)
 =2\operatorname{tr}\bigl((BB^*)^2\bigr)
 \le 2\norm B_{\op}^2\norm B_{\HS}^2.
\end{equation}
Equations~\eqref{eq:local-risk},~\eqref{eq:op}, and~\eqref{eq:variance} imply that the variance of $h^{-p}\norm{\phi^{1/2}X_h}^2$ is $O(h^{2-p})$. This bound is summable along the dyadic meshes. Chebyshev's inequality and Borel--Cantelli yield the limit against each fixed nonnegative $\phi$. A countable dense set of continuous tests and the bound on the total mass yield~\eqref{eq:weak}.

Write $e_m=u_{h_m}-Y_m$ and define
\begin{equation*}
 Z_m=
 \begin{cases}
  \dfrac{\ip{X_{h_m}}{e_m}}{h_m^{p/2}\norm{e_m}}, & e_m\ne0,\\
  0, & e_m=0.
 \end{cases}
\end{equation*}
Given the observations and any independent auxiliary randomness, the real inner product $\ip{X_{h_m}}{e_m}$ is centered Gaussian with variance at most $Ch_m^2\norm{e_m}^2$. Hence $Z_m$ is conditionally centered and has conditional variance $O(h_m^{2-p})$. The same summability argument gives $Z_m\to0$ almost surely.

Expanding $\norm{X_{h_m}+e_m}^2$ and using $2h_m^{p/2}\norm{e_m}\le h_m^p+\norm{e_m}^2$ gives~\eqref{eq:path-pythagoras}. The optimality assertions follow by dividing by $h_m^p$.
\end{proof}

\section{The cell-averaged FFT method}\label{sec:solver}

We construct a numerical solution of the approximate problem. The double cell average of the Green function depends only on the displacement between the source and output cells. We therefore use the FFT to evaluate the discrete convolution. We estimate the error due to the $L^2$ projection onto the output mesh separately from the errors due to the approximation of the source coefficient and the convolution weights.

Throughout this section, $D$ and $Q$ are aligned boxes. The source mesh size is $h$, and the output mesh size is $H=h/r$ for a positive integer $r$. Each box is a union of cubes from its respective mesh. The meshes have a common origin, so each coarse cube is a union of fine cubes. Let $P_H^Q$ be the $L^2(Q)$ orthogonal projection onto functions that are constant on each output cell. We set $V_{h,H}=P_H^Q u_h$.

Let $a_H$ be constant on the fine source cells and satisfy $\norm{a_H-a}_{L^\infty(D)}\le CH$. Both cell averages and values at the cell centers satisfy this estimate. On each fine source cell, the discretized forcing is $a_H$ times the value of $P_hW$ on the coarse cell containing it. With exact convolution weights, the numerical method gives $P_H^Q T M_{a_H}P_hW$. The numerical solution $U_{h,H}$ uses the approximate weights of Section~\ref{sec:weight-construction}. On a sequence of meshes, write $U_m=U_{h_m,H_m}$ and $V_m=V_{h_m,H_m}$.

\subsection{The output projection}

The difference $u_h-V_{h,H}$ is the error due to the $L^2$ projection onto the output mesh. Let $J=K*\mathbf 1_{C_0}$. Let $\Pi_{1/r}$ be averaging on the tessellation of $\R^d$ by cubes of side $1/r$ whose boundaries align with those of $C_0$. Define
\begin{equation}\label{eq:penalty}
 B_r=\norm{(I-\Pi_{1/r})J}_{L^2(\R^d)}^2 .
\end{equation}

\begin{lemma}\label{lem:penalty}
The constants $B_r$ are positive and $r^2 B_r\to \norm{\nabla J}_{L^2(\R^d)}^2/12$ as $r\to\infty$. For each $\phi\in C(\overline Q)$ with $\phi\ge0$, uniformly over positive integers $r$ and $H=h/r$,
\begin{equation}\label{eq:output-limit}
 h^{-p}\norm{\phi^{1/2}(I-P_H^Q)T_\mu P_h}_{\HS}^2
 =B_r\int_D\phi(x)\mu(x)\,dx+o(r^{-2})
 \quad\text{as }h\to0 .
\end{equation}
\end{lemma}

\begin{proof}
The Fourier transform of $\nabla J$ has absolute value
\begin{equation*}
 \abs{\xi}^{1-2\alpha}\prod_{\ell=1}^d\abs{\sinc(\xi_\ell/2)} .
\end{equation*}
It is square integrable at the origin by $p<2$. At infinity, $\abs{\xi}^{1-2\alpha}$ is bounded and the sinc product is in $L^2$. Poincar\'e's inequality gives $B_r\le Cr^{-2}$. If $B_r=0$, the locally $H^1$ function $J$ is constant on every cell and its weak gradient vanishes. This contradicts positivity and decay of $J$.

Put $\delta=1/r$. On each cell of width $\delta$ centered at $z$, approximate $(I-\Pi_\delta)J$ by $(\Pi_\delta\nabla J)\cdot(x-z)$. The difference has zero cell mean and gradient $(I-\Pi_\delta)\nabla J$. Its norm over all cells is at most $C\delta\norm{(I-\Pi_\delta)\nabla J}_{L^2}=o(\delta)$. The squared norm of the linear term is $\delta^2\norm{\Pi_\delta\nabla J}_{L^2}^2/12$. This proves the limit of $r^2 B_r$.

We estimate the error due to replacing $a$ by $a_h$ uniformly in $r$. Lemma~\ref{lem:riesz} and output-cell Poincar\'e give
\begin{equation*}
 \norm{(I-P_H^Q)T M_{a-a_h}P_h}_{\HS}
 \le CH\norm{M_{a-a_h}P_h}_{\HS}
 \le CHh^{1-d/2}
 =o(h^{p/2}/r).
\end{equation*}
In the second inequality we sum over the $\abs D h^{-d}$ normalized cell indicators and use $\norm{a-a_h}_\infty\le Ch$. The last step uses $\alpha<1$. Replacing $G$ by $K$ in the operator with coefficient $a_h$ costs $O(H)$ in the same norm by the local $H^1$ bound on $G_{\rm rem}$ and output-cell Poincar\'e. Multiplication by $\phi^{1/2}$ preserves these orders.

On a normalized indicator of a source cell $E$ centered at $c$, the frozen Riesz operator produces $a_E h^{2\alpha-d/2}J((x-c)/h)$, where $a_E$ is the value of $a_h$ on $E$. Its contribution to the squared weighted projection error is
\begin{equation*}
 a_E^2 h^{p+d}\int_{\R^d}\widetilde\phi(c+ht)
      \abs{(I-\Pi_{1/r})J(t)}^2\,dt .
\end{equation*}
Here $\widetilde\phi$ is the zero extension outside $Q$. On bounded sets of $t$, we may replace $\widetilde\phi(c+ht)$ by $\phi(c)$ with an error tending uniformly to zero for $c\in D$. For $L>\sqrt d$, the integral of $\abs{(I-\Pi_{1/r})J}^2$ over $\abs t>L$ is bounded by $Cr^{-2}\int_{\abs t>L-\sqrt d}\abs{\nabla J(t)}^2dt$ by Poincar\'e on the cells meeting the tail. Thus the replacement has error $o(r^{-2})$ uniformly in $r$ and $c$. Summing the cell contributions gives $h^p B_r\int_D\phi\mu+o(h^p/r^2)$ by Riemann sums. The frozen Riesz norm is $O(h^{p/2}/r)$, so the estimates for the replacements of the source coefficient and the Green function also give $o(h^p/r^2)$ in the squared norm. This proves~\eqref{eq:output-limit}.
\end{proof}

\subsection{Construction and cost of the weights}\label{sec:weight-construction}

For a fine-grid offset $\mathbf m\in\mathbb Z^d$, average $G$ over two cubes of side $H$ whose centers differ by $H\mathbf m$. The resulting weight is
\begin{equation}\label{eq:tent}
 \mathcal A_H(\mathbf m)=\int_{(-1,1)^d}G\bigl(H(\mathbf m+z)\bigr)\prod_{\ell=1}^d(1-\abs{z_\ell})\,dz .
\end{equation}
The tent factor is the density of the difference of two independent uniform points in $C_0$. Let $\widetilde{\mathcal A}_H$ be the approximation to $\mathcal A_H$ constructed below. We index the output and fine source cells by $i,j\in\mathbb Z^d$, respectively. If $f^H_j$ is the source density on cell $j$, the computed value on cell $i$ is
\begin{equation*}
 H^d\sum_j\widetilde{\mathcal A}_H(i-j)f^H_j .
\end{equation*}
We denote this convolution operator by $\widetilde T_H$. In the present model the source strength enters through $f^H_j$, so $U_{h,H}=\widetilde T_H M_{a_H}P_hW$. The weights depend only on the grid and the free-space Green function. Zero padding gives the linear convolution without wraparound. For free-space FFT convolution with truncated Green functions, see~\cite{VGF}. Let $N$ denote the number of padded points. The set of offsets has $O(N)$ elements and $N\asymp H^{-d}$ for the fixed boxes.

The fractional-power resolvent identity~\cite[Eq.~(5.28)]{MS} separates the classical outgoing Helmholtz term from a Stieltjes integral. See also~\cite[Sections~2.3 and~2.5]{ZCV} for the free-space outgoing Green function in two and three dimensions. Let $g_k$ be the classical outgoing Helmholtz kernel and $\mathcal Y_t$ the kernel of $(-\Delta+t)^{-1}$ for $t>0$. Put
\begin{equation}\label{eq:stieltjes}
 c_k=\frac{k^{2-2\alpha}}{\alpha},\qquad
 b(t)=\frac{\sin(\pi\alpha)}{\pi}\,
 \frac{t^\alpha}{t^{2\alpha}-2k^{2\alpha}t^\alpha\cos(\pi\alpha)+k^{4\alpha}} .
\end{equation}
Then
\begin{equation}\label{eq:heat-split}
 G(x)=c_k\bigl(g_k(x)-\mathcal Y_1(x)\bigr)
    +\int_0^\infty \omega(s) (4\pi s)^{-d/2}e^{-\abs x^2/(4s)}\,ds ,
 \qquad
 \omega(s)=c_k e^{-s}+\int_0^\infty b(t)e^{-st}\,dt .
\end{equation}
Indeed, the function $(z^\alpha-k^{2\alpha})^{-1}-c_k(z-k^2)^{-1}$ is analytic off the nonpositive real axis after the singularity at $k^2$ is removed. It decays at infinity. Cauchy's formula about the cut represents it as $\int_0^\infty b(t)/(z+t)\,dt$ using the jump in its boundary values. Taking the outgoing boundary value at $z=\abs{\xi}^2$ gives the corresponding free resolvent decomposition. The heat representation of $\mathcal Y_t$ and Tonelli's theorem for its positive mixture give~\eqref{eq:heat-split} as an identity of locally integrable kernels. The term $g_k$ carries the outgoing prescription. The three-dimensional radial form of the decomposition is given in~\cite[Eqs.~(2.8)--(2.9)]{QXZfree} and~\cite[Lemma~2.1]{QXZinhom}. Related Gaussian representations of the classical Helmholtz Green function appear in work on fast convolution~\cite{BKM} and on the separation of oscillatory and singular components~\cite{Beylkin}.

For positive real $s$, define the one-dimensional average
\begin{equation}\label{eq:heat-average}
\begin{aligned}
 q_H(j,s)&=\frac{F_s((j+1)H)-2F_s(jH)+F_s((j-1)H)}{H^2},\\
 F_s(x)&=\frac{x}{2}\operatorname{erf}\bigl(x/(2\sqrt s)\bigr)
       +\sqrt{s/\pi}\,e^{-x^2/(4s)}.
\end{aligned}
\end{equation}
Here $F_s''$ is the one-dimensional heat kernel. Consequently, the heat part of $\mathcal A_H(\mathbf m)$ is the integral of $\omega(s)\prod_\ell q_H(\mathbf m_\ell,s)$. At each quadrature node the weights are formed as a product of one-dimensional arrays.

\begin{theorem}\label{thm:weights}
In exact arithmetic, there are deterministic choices of quadrature for~\eqref{eq:heat-split} that give weights $\widetilde{\mathcal A}_H(\mathbf m)$ satisfying
\begin{equation}\label{eq:uniform-weights}
 \sup_{\mathbf m}\abs{\widetilde{\mathcal A}_H(\mathbf m)-\mathcal A_H(\mathbf m)}\le CH
\end{equation}
on the prescribed displacement set. With $\Lambda_H=1+\log(1/H)$ for $H\le 1$, the construction uses $O(N\Lambda_H^2+\Lambda_H^4)$ scalar operations and special-function evaluations. The convolution, including its FFT preparation, uses an additional $O(N\log N)$ operations.
\end{theorem}

\begin{proof}
The radial function $c_k(g_k-\mathcal Y_1)$ is Lipschitz on bounded intervals in the radial variable $\rho$. In dimension three this follows from $(e^{ik\rho}-e^{-\rho})/(4\pi\rho)$ and its continuous limit at zero. In dimension two the common logarithmic term of $iH_0^{(1)}(k\rho)/4$ and $K_0(\rho)/(2\pi)$ cancels. The remainder is a smooth function of $\rho^2$ plus a term of the form $\rho^2\log\rho$ times a smooth function of $\rho^2$~\cite[Sections~10.8 and~10.31]{DLMF}. Evaluating the remainder at the displacement center thus introduces error at most $L_G H\sqrt{d/6}$, where $L_G$ is a Lipschitz bound. This bound uses the second moment of the tent density in~\eqref{eq:tent}. The two-dimensional power series give uniform error $O(H)$ on bounded displacement sets with $O(\Lambda_H)$ terms, by the factorial decay of their coefficients.

Fix $\theta\in(0,\pi/2)$. In $\abs{\arg t}\le\theta$ the denominator factors of $b(t)$ have magnitude bounded below by a constant times $\abs t^\alpha+k^{2\alpha}$. Thus
\begin{equation*}
 \abs{b(t)}\le C\min(\abs t^\alpha,\abs t^{-\alpha}).
\end{equation*}
For $\abs{\arg s}\le\theta$, this gives $\abs{\omega(s)}\le C\min(\abs s^{\alpha-1},\abs s^{-1-\alpha})$. On the same sector the modulus of the heat kernel is bounded by a constant times the heat kernel at positive real time $\abs s/\cos\theta$. Averaging gives
\begin{equation*}
 \left|\prod_{\ell=1}^d q_H(\mathbf m_\ell,s)\right|
 \le C\min(H^{-d},\abs s^{-d/2}).
\end{equation*}
After $s=e^v$, the integrand is analytic in the strip $\abs{\operatorname{Im}v}<\pi/2$ and its $L^1$ norms on the lines $\operatorname{Im}v=\pm\theta$ are bounded by $C H^{2\alpha-d}$. This follows by splitting the integral at $\abs s=H^2$ and using the $s^{\alpha-1}$ estimate for $\omega$. The infinite trapezoidal rule with step $\delta$ has error at most $C H^{2\alpha-d}/(\exp(2\pi\theta/\delta)-1)$ by the strip estimate~\cite[Theorem 5.1]{TW}. The discarded left and right sums outside $(-\ell_H,\ell_H)$ are bounded by $CH^{-d}e^{-\alpha\ell_H}$ and $Ce^{-(\alpha+d/2)\ell_H}$. Choosing $\ell_H=c_1 \Lambda_H$ and $\delta=c_2/\Lambda_H$ with suitable constants gives error $O(H)$ and $O(\Lambda_H^2)$ nodes.

For the inner integral in $\omega(s)$, use $t=e^y/s$ at each positive $s$. For related logarithmic quadrature of resolvent integrals for fractional powers, see~\cite{BP}. After division by $s^{\alpha-1}$, the integrand is bounded on $\abs{\operatorname{Im}y}\le\theta$ by
\begin{equation*}
 C \exp\bigl((1-\alpha)\operatorname{Re}y-\cos(\theta)e^{\operatorname{Re}y}\bigr).
\end{equation*}
The same trapezoidal estimate with a suitable step of order $\Lambda_H^{-1}$ yields error bounded by $H^{1+d-2\alpha}s^{\alpha-1}$ when the range is a sufficiently large multiple of $\Lambda_H$. A common uniform grid in $\log t$ is equally valid, since the strip estimate is translation invariant and $y=\log t+\log s$ shifts the log nodes by at most $O(\Lambda_H)$ on the retained outer range. In the outer rule the discrete sum of $s^\alpha\min(H^{-d},s^{-d/2})$ times the log step is $O(H^{2\alpha-d})$, by splitting the geometric sums at $H^2$. The inner errors therefore contribute $O(H)$. The two rules use $O(\Lambda_H^4)$ operations to obtain the heat coefficients. At each outer node, the arrays in~\eqref{eq:heat-average} and their products can be formed for all displacements in $O(N)$ operations. This proves the weight error and the claimed work bounds.
\end{proof}

\subsection{Almost-sure error estimates for the FFT method}

We estimate the error in exact arithmetic. The error $U_{h,H}-V_{h,H}$ has two contributions, one from approximating the source coefficient and one from approximating the convolution weights. The operator $P_H^Q T M_{a_H-a}P_h$ has Hilbert--Schmidt norm at most $\norm T_{\HS}\norm{a_H-a}_{L^\infty(D)}$. Let $\eta$ be the maximum of $\abs{\widetilde{\mathcal A}_H(i-j)-\mathcal A_H(i-j)}$ over all output indices $i$ and fine source indices $j$. The operator associated with the weight error has Hilbert--Schmidt norm at most $\eta |Q|^{1/2}\norm{a_H}_{L^2(D)}$. Consequently,
\begin{equation}\label{eq:weight-error}
 \E\norm{U_{h,H}-V_{h,H}}^2\le C(H^2+\eta^2) .
\end{equation}
Theorem~\ref{thm:weights} provides $\eta=O(H)$. The output projection error in Lemma~\ref{lem:penalty} therefore determines the leading term of $\norm{u_h-U_{h,H}}$ along dyadic refinements.

Let $S_H$ be the range of $P_H^Q$ and set $C_\mu(r)=C_\mu+I_\mu B_r$. The following theorem applies the almost-sure error expansion of Section~\ref{sec:sharp} to the choice $Y_m=U_m$ and to other approximations in the same output space.

\begin{theorem}\label{thm:algorithm}
Let $h_m=h_0 2^{-m}$ and $H_m=h_m/r_m$ with deterministic positive integers $r_m$. Let $\eta_m$ denote the maximum absolute error in the weights over all required offsets at level $m$. Suppose that $\eta_m=O(H_m)$. Then almost surely,
\begin{equation}\label{eq:fixed-r}
 \norm{u_{h_m}-U_m}^2\sim I_\mu h_m^p B_{r_m},\qquad
 h_m^{-p}\norm{u-U_m}^2=C_\mu(r_m)+o(1).
\end{equation}
For any fixed sequence of $\mathcal F_{h_m}$-measurable approximations $Y_m\in S_{H_m}$, one has almost surely
\begin{equation}\label{eq:relative-error}
 \frac{\norm{u-Y_m}^2}
 {C_\mu(r_m)h_m^p+\norm{Y_m-U_m}^2}
 \longrightarrow 1 .
\end{equation}
In particular, $\liminf_m(h_m^{-p}\norm{u-Y_m}^2-I_\mu B_{r_m})\ge C_\mu$. Convergence of $h_m^{-p}\norm{u-Y_m}^2$ to $C_\mu$ along the full sequence requires $r_m\to\infty$. Under this condition, $h_m^{-p}\norm{u-U_m}^2\to C_\mu$ almost surely.
\end{theorem}

\begin{proof}
By~\eqref{eq:weight-error} and summability of $h_m^{2-p}$, Markov's inequality and the Borel--Cantelli lemma give $\norm{U_m-V_m}=o(h_m^{p/2}/r_m)$ almost surely. The coefficient $(I-P_H^Q)T_\mu P_h$ has operator norm $O(H)$ by Lemma~\ref{lem:riesz} and Poincar\'e. Its squared Hilbert--Schmidt norm is $O(h^p/r^2)$ by Lemma~\ref{lem:penalty}. Thus the variance of the squared residual normalized by $h^p/r^2$ is $O(h^{2-p})$ by~\eqref{eq:variance}. It follows from~\eqref{eq:output-limit} and $B_r\asymp r^{-2}$ that $\norm{u_{h_m}-V_m}^2\sim I_\mu h_m^p B_{r_m}$ almost surely. The estimate on $\norm{U_m-V_m}$ gives the first assertion in~\eqref{eq:fixed-r}.

For $Y_m\in S_{H_m}$, spatial orthogonality gives
\begin{equation*}
 \norm{u_{h_m}-Y_m}^2
 =\norm{u_{h_m}-V_m}^2+\norm{V_m-Y_m}^2 .
\end{equation*}
Since $\norm{U_m-V_m}=o(h_m^{p/2})$, replacing $V_m$ by $U_m$ in the second term changes it by $o(h_m^p+\norm{U_m-Y_m}^2)$. The expansion~\eqref{eq:path-pythagoras} now yields~\eqref{eq:relative-error}, since $C_\mu(r_m)$ is bounded above and below by positive constants. Taking $Y_m=U_m$ yields the second assertion of~\eqref{eq:fixed-r}. Dropping the nonnegative term in the denominator of~\eqref{eq:relative-error} gives the lower bound.

If $r_m$ has a bounded subsequence, then $B_{r_m}$ is bounded away from zero on that subsequence and the normalized error of any fixed $Y_m$ stays above $C_\mu$ there by a positive amount for sufficiently large $m$. If $r_m\to\infty$, then $B_{r_m}\to0$ and~\eqref{eq:fixed-r} gives the full-sequence limit for $U_m$.
\end{proof}

If the ratios are unbounded and have a bounded subsequence, then~\eqref{eq:fixed-r} gives
\begin{equation*}
 \liminf_{m\to\infty} h_m^{-p}\norm{u-U_m}^2=C_\mu .
\end{equation*}
The normalized error stays above $C_\mu$ by a positive amount along the bounded subsequence.

A piecewise constant output on the source mesh incurs an extra leading term $I_\mu B_1$. With $n=\abs D h^{-d}$ white noise integrals on the source mesh and a bounded-cost evaluation of the source coefficient per fine source cell, the work including weight generation is $O(nr^d(1+\log(nr^d))^2)$. For fixed $r$, the leading constant is $C_\mu(r)$, and the excess $C_\mu(r)-C_\mu$ is of order $r^{-2}$. The numerical solution attains the optimal constant $C_\mu$ along any deterministic sequence $r_m\to\infty$. For instance, taking $r_m$ dyadic and of order $1+\log(m+1)$ gives work $O(n(1+\log n)^2(1+\log\log n)^d)$ as $n\to\infty$.

\subsection{Differences between successive approximations}\label{sec:adjacent}

To compare the numerical results with the error estimates, we use differences between successive approximations of the same realization. On nested output meshes, the numerical solution on the coarse mesh is an admissible approximation on the fine mesh. Applying~\eqref{eq:relative-error} to this choice gives the following limit.

\begin{corollary}\label{cor:adjacent}
Under the assumptions of Theorem~\ref{thm:algorithm}, if $S_{H_m}\subset S_{H_{m+1}}$ for all $m$, then almost surely
\begin{equation}\label{eq:adjacent}
 \frac{\norm{U_{m+1}-U_m}^2}{(1-2^{-p})h_m^p}
 =
 \frac{C_\mu(r_m)-2^{-p}C_\mu(r_{m+1})}{1-2^{-p}}+o(1).
\end{equation}
At fixed $r_m=r$, the ratio tends to $C_\mu(r)$. The ratio for the conditional means $u_{h_m}$ and $u_{h_{m+1}}$ tends to $C_\mu$ almost surely. The same limit holds for the FFT approximations when $r_m\to\infty$, with arbitrary admissible output grids at each level.
\end{corollary}

\begin{proof}
The coarse field $U_m$ is $\mathcal F_{h_{m+1}}$-measurable because the source data are nested. The inclusion of output spaces gives $U_m\in S_{H_{m+1}}$. At level $m+1$, choose $Y_{m+1}=U_m$ in~\eqref{eq:relative-error}. The squared error on the coarse mesh in~\eqref{eq:fixed-r} is $O(h_m^p)$ almost surely, so
\begin{equation*}
 C_\mu(r_{m+1})h_{m+1}^p+\norm{U_{m+1}-U_m}^2
 =\norm{u-U_m}^2+o(h_m^p).
\end{equation*}
Substituting the error on the coarse mesh from~\eqref{eq:fixed-r} gives~\eqref{eq:adjacent}. For the conditional means the same argument uses $Y_{m+1}=u_{h_m}$ in~\eqref{eq:path-pythagoras}. If $r_m\to\infty$, then~\eqref{eq:fixed-r} gives $\norm{U_m-u_{h_m}}=o(h_m^{p/2})$ almost surely. The difference between successive numerical solutions therefore has the same leading squared norm as the difference of the conditional means.
\end{proof}

For the experiment on a common output mesh, we compute the two numerical solutions at source mesh sizes $h$ and $h/2$ on a common output mesh of size $H=h/(2r)$. Let
\begin{equation*}
 L_{h,H}=P_H^Q T_\mu(P_{h/2}-P_h)W,\qquad
 \widetilde L_{h,H}=U_{h/2,H}-U_{h,H}.
\end{equation*}
Both numerical solutions use the same $a_H$ and the same convolution weights. Thus $\widetilde L_{h,H}$ can be computed by applying the convolution to $a_H(P_{h/2}-P_h)W$. Orthogonality of the source and output projections gives
\begin{equation}\label{eq:common-trace}
 \E\norm{L_{h,H}}^2=R_h-R_{h/2}
  -\norm{(I-P_H^Q)T_\mu P_{h/2}}_{\HS}^2
  +\norm{(I-P_H^Q)T_\mu P_h}_{\HS}^2 .
\end{equation}
The projection terms have refinement ratios $r$ and $2r$, respectively. For fixed $r$, dividing $\E\norm{L_{h,H}}^2$ by $(1-2^{-p})h^p$ gives the limit
\begin{equation}\label{eq:common-limit}
 I_\mu\kappa_r,\qquad
 \kappa_r=\mathcal C_{d,\alpha}-\frac{2^{-p}B_r-B_{2r}}{1-2^{-p}}.
\end{equation}
The output projection reduces the limiting constant by $C_\mu-I_\mu\kappa_r$. This reduction is at most $I_\mu 2^{-p}B_r/(1-2^{-p})$. The cube Poincar\'e inequality gives $B_r\le \norm{\nabla J}_{L^2(\R^d)}^2/(\pi r)^2$, so the limiting constant tends to $C_\mu$ as $r\to\infty$.

The squared difference on the common output mesh also has a weak limit after normalization. Along the dyadic grids, for fixed $r$ one has almost surely
\begin{equation}\label{eq:common-local}
 \frac{|\widetilde L_{h,H}|^2}{(1-2^{-p})h^p}\,dx
 \ \rightharpoonup\ \kappa_r\mu\mathbf 1_D\,dx .
\end{equation}
To prove this, take $\phi\in C(\overline Q)$ with $\phi\ge0$ and put $\phi_H=P_H^Q\phi$. The integral of $|L_{h,H}|^2$ against $\phi$ equals that against $\phi_H$. Since $\phi_H$ is constant on output cells, the weighted mean satisfies~\eqref{eq:common-trace} with $\phi_H^{1/2}$ multiplying each coefficient, including those in $R_h$ and $R_{h/2}$. We can replace $\phi_H$ by $\phi$ in these traces with error $o(h^p)$ by uniform continuity and the $O(h^p)$ bounds for all four unweighted terms. Equations~\eqref{eq:local-risk} and~\eqref{eq:output-limit} give the limiting weighted mean. The operator defining $L_{h,H}$ has norm $O(h)$ by~\eqref{eq:op}. By~\eqref{eq:variance}, the normalized squared norm has variance $O(h^{2-p})$. Summability and a countable dense set of continuous test functions give the almost-sure weak limit for $L_{h,H}$. The estimate in~\eqref{eq:weight-error} applies with $P_{h/2}-P_h$ in place of $P_h$, using the same $a_H$ for the two levels. The resulting estimate $\norm{\widetilde L_{h,H}-L_{h,H}}=o(h^{p/2})$ holds almost surely. The Cauchy--Schwarz inequality then gives the same weak limit for the normalized squared difference computed by the FFT. For a deterministic sequence $r_m\to\infty$, the same proof with $h=h_m$ and $H_m=h_m/(2r_m)$ gives the almost-sure weak limit $\mathcal C_{d,\alpha}\mu\mathbf 1_D\,dx$, by the uniform estimate in~\eqref{eq:output-limit}.

For each positive integer $r$, one has $0<\kappa_r<\mathcal C_{d,\alpha}$. Take an orthonormal basis of the mean-zero functions on $C_0$ that are constant on its half-size children, extended by zero. For each basis function $q$, the Fourier transform of $K*q$ is $|\xi|^{-2\alpha}\widehat q(\xi)$. The cancellation $\widehat q(\xi)=O(|\xi|)$ at zero, together with $p<2$, gives square integrability near zero. At high frequencies, $2\alpha>1$ gives $K*q\in H^1(\R^d)$.

For $\mu=1$, replacing $G$ by $K$ in the operator defining $L_{h,H}$ changes it by $O(h)$ in Hilbert--Schmidt norm. This follows from the smooth remainder estimate in the proof of Theorem~\ref{thm:risk} and the identity $P_{h/2}-P_h=(I-P_h)(P_{h/2}-P_h)$. Rescale each parent cell to $C_0$. The rescaled observation regions contain balls of radii tending uniformly to infinity because $\overline D\subset Q$. Summing the projected and unprojected energies over the rescaled basis of mean-zero piecewise constant functions gives
\begin{align*}
 (1-2^{-p})\kappa_r&=\sum_q\norm{\Pi_{1/(2r)}(K*q)}_{L^2(\R^d)}^2,\\
 (1-2^{-p})\mathcal C_{d,\alpha}&=\sum_q\norm{K*q}_{L^2(\R^d)}^2 .
\end{align*}
Pairing $\Pi_{1/(2r)}(K*q)$ with $q$ gives the positive integral $(2\pi)^{-d}\int_{\R^d}|\xi|^{-2\alpha}|\widehat q(\xi)|^2\,d\xi$, since $q$ belongs to the range of the projection. Also, a function in $H^1(\R^d)$ that is constant on every grid cell has weak gradient zero and must vanish. Hence each nonzero response $K*q$ has a nonzero projection residual. The two identities give the strict bounds on $\kappa_r$. In particular, for a subregion of $D$ whose boundary has measure zero, the fraction of $\norm{\widetilde L_{h,H}}_{L^2(D)}^2$ contributed by the subregion tends to the corresponding fraction of the integrated source strength.

A known deterministic forcing $f\in L^2(D)$ can be included in~\eqref{eq:model}. The solution and its conditional mean both contain $Tf$, so the centered residual is $T_\mu(I-P_h)W$. In the FFT method, use the cell averages of $f$ on the output mesh. Input projection, output projection, and the weight bound together perturb $Tf$ by $O(H)\norm f_{L^2(D)}$, which is of smaller order than $h^{p/2}/r$.

\section{Numerical experiments}\label{sec:experiment}

We present numerical examples in two and three dimensions. In each dimension, we generate the white noise integrals on the finest source mesh and sum them to obtain the integrals on each coarser mesh. We compare the differences between successive numerical solutions from this single realization with the limits in Section~\ref{sec:adjacent}. This comparison requires no evaluation of the exact solution.

\subsection{Constant source strength}

For the tests with constant source strength, we set $\mu=1$. We take $D=(0,1)^d$, $Q=(-1/4,5/4)^d$, and $k=1$. The pairs $(d,\alpha)=(2,3/4)$ and $(3,15/16)$ give $p=1$ and $3/4$, respectively. On the finest source mesh, the cell integrals are independent normal random variables with mean zero and variance equal to the cell volume. This mesh has $256$ cells per coordinate direction in two dimensions and $64$ in three dimensions. We use the PCG64 generator~\cite{ONeill} with the dimension and the master seed $20261119$.

For source mesh sizes $h$ and $h/2$, we compare two choices of output meshes. With a separate output mesh for each approximation, the mesh sizes are $h/r$ and $h/(2r)$. The coarse output is constant on each of its fine subcells, so we can take the difference on the finer grid. Dividing its squared norm by $(1-2^{-p})h^p$ gives the limit $C_\mu(r)$ at fixed $r$ by Corollary~\ref{cor:adjacent}. The additional term measures the effect of representing each field on its own output grid.

When both levels use the common finer output grid $H=h/(2r)$, linearity permits a single convolution with the difference of the two source densities. Output averaging then removes part of the energy of this difference, as quantified by~\eqref{eq:common-trace} and~\eqref{eq:common-limit}. With $\widetilde L_{h,H}$ denoting the computed difference on this grid, we measure
\begin{equation}\label{eq:measured}
 S_{h,r}=\frac{\norm{\widetilde L_{h,H}}^2}{(1-2^{-p})h^p}.
\end{equation}
The leading constants $C_\mu$ for the conditional mean are $0.0829842520$ in two dimensions and $0.0387168158$ in three dimensions, from~\eqref{eq:constant-cell}. We evaluate the fixed-$r$ limits for both choices by quadrature of the Riesz profile formula in Appendix~\ref{app:quadrature}. We also compute $\E S_{h,r}$ directly on the finite grid. For a parent source cell of side $h$, let $E_j$ denote its children of side $h/2$. Let $v_j$ be the numerical solution on the common output grid for the normalized input $|E_j|^{-1/2}\mathbf 1_{E_j}$, and let $\overline v$ be the mean of these solutions. The parent cell contributes $\sum_j\norm{v_j-\overline v}^2$ to the covariance trace.

\begin{figure}[!htbp]
 \centering
 \includegraphics[width=\linewidth]{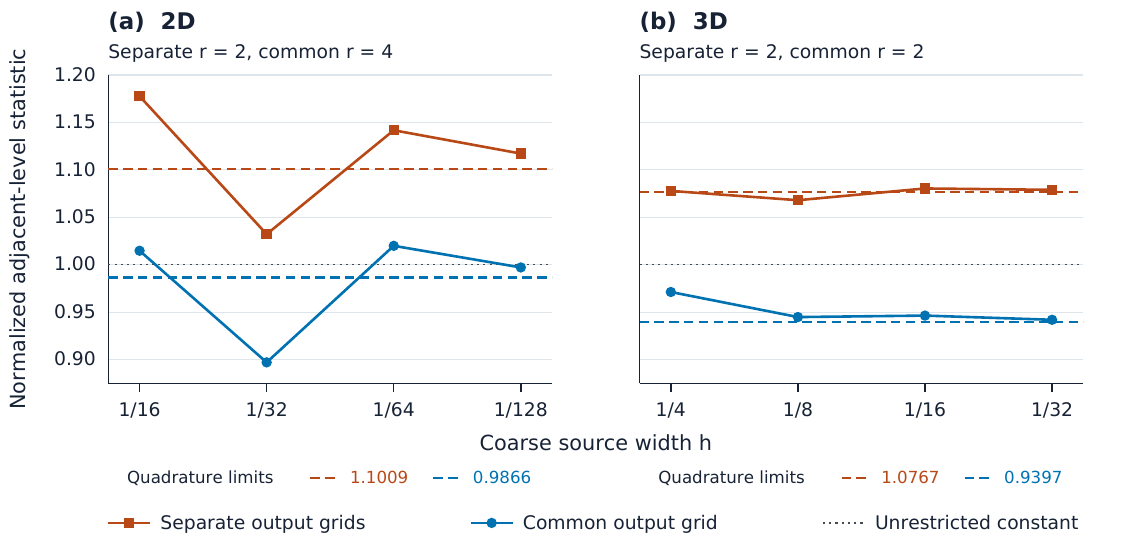}
 \caption{Squared $L^2$ norms of the differences between successive numerical solutions, divided by $(1-2^{-p})h^pC_\mu$. The two output mesh choices use the same white noise realization in each dimension. Each panel shows the four finest source mesh sizes shared by the two tests. In two dimensions the refinement ratios are $r=2$ for separate output grids and $r=4$ for the common grid. In three dimensions both ratios are $r=2$. The dotted line marks $1$, and the dashed lines give quadrature approximations to the normalized fixed-ratio limits.}
 \label{fig:single}
\end{figure}

Figure~\ref{fig:single} compares the two output grid choices. The fixed-ratio limit lies above $1$ when each level uses its own output grid and below $1$ when the two levels share the fine output grid. Table~\ref{tab:accuracy} gives the finite-grid means and fixed-ratio limits for the common-grid experiment. In three dimensions, increasing $r$ from $2$ to $4$ at $h=1/32$ raises the measured ratio by $0.04288$ (predicted $0.04271$). Increasing $r$ from $4$ to $8$ at $h=1/16$ gives an increment of $0.01324$ (predicted $0.01299$). The finest-grid fractions of energy outside $D$ are $1.63\%$ in two dimensions at $r=8$ and $4.47\%$ in three dimensions at $r=4$.

\begin{table}[!htbp]
\centering
\caption{Common-grid statistics, their finite-grid means, and quadrature approximations to the limits in~\eqref{eq:common-limit}, all normalized by $C_\mu$.}
\label{tab:accuracy}
\small
\begin{tabular}{@{}rrrrrr@{}}
\toprule
$d$ & $h$ & $r$ & $S_{h,r}/C_\mu$ & $\E S_{h,r}/C_\mu$ & Limit (quadrature) \\
\midrule
2 & $1/128$ & 4 & 0.99717 & 0.99081 & 0.9866 \\
2 & $1/128$ & 8 & 1.00692 & 1.00067 & 0.9964 \\
3 & $1/16$ & 2 & 0.94641 & 0.93942 & 0.9397 \\
3 & $1/16$ & 4 & 0.98997 & 0.98213 & 0.9824 \\
3 & $1/16$ & 8 & 1.00321 & 0.99512 & 0.9954 \\
3 & $1/32$ & 2 & 0.94198 & 0.93983 & 0.9397 \\
3 & $1/32$ & 4 & 0.98485 & 0.98254 & 0.9824 \\
\bottomrule
\end{tabular}
\end{table}

Adding $C_\mu-I_\mu\kappa_r$ to $S_{h,r}$ restores the limit $C_\mu$ for fixed $r$. At the finest $h$ in each dimension, the corrected ratios are $1.01061$ and $1.01051$ in two dimensions at $r=4,8$, and $1.00228$ and $1.00244$ in three dimensions at $r=2,4$.

For the experiment with each level on its own output grid, the finest-grid ratios at $r=1$ and $2$ are $1.3372$ and $1.1172$ in two dimensions. The quadrature approximations to $C_\mu(r)/C_\mu$ are $1.3054$ and $1.1009$, respectively. In three dimensions the measured ratios are $1.2286$ and $1.0790$, with predictions $1.2271$ and $1.0767$. In both dimensions the results illustrate the additional error due to the output projection in Theorem~\ref{thm:algorithm}.

\subsection{Spatially varying source strength}

In the two-dimensional experiment with $\alpha=3/4$ and $k=1$, the source strength functions are
\begin{equation*}
 \mu_0(x)=1,\qquad \mu_1(x)=1+\tfrac12\sin(2\pi x_1),\qquad x\in D.
\end{equation*}
Both have integral $1$, so they share the same $C_\mu$ and the same fixed-$r$ limit in~\eqref{eq:common-limit}.

\begin{figure}[!htbp]
 \centering
 \includegraphics[width=\linewidth]{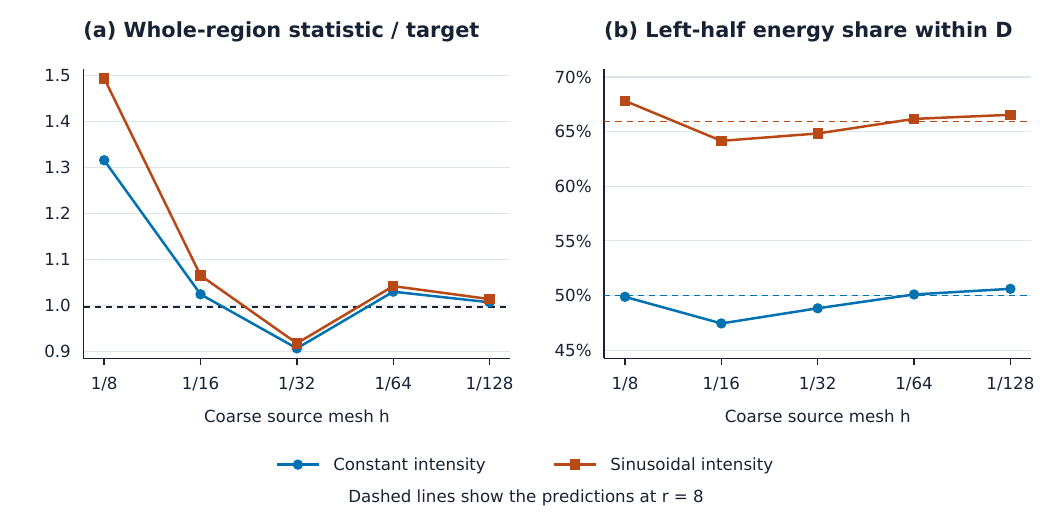}
 \caption{Results on a common output mesh with $r=8$. The left panel shows $S_{h,r}/C_\mu$ for two source strengths and the same white noise realization. The right panel shows the fraction of the squared $L^2(D)$ norm of the difference between successive numerical solutions contributed by the left half of $D$. Dashed lines indicate the limiting values.}
 \label{fig:intensity}
\end{figure}

The two numerical examples use the same white noise integrals and the same convolution weights. For $\mu_1$, the cell averages of the coefficient $\sqrt{\mu_1}$ are evaluated by a 12-point Gauss--Legendre rule in the $x_1$ direction.

Let $D_L=(0,1/2)\times(0,1)$. The right panel of Figure~\ref{fig:intensity} gives the ratio of $\norm{\widetilde L_{h,H}}_{L^2(D_L)}^2$ to $\norm{\widetilde L_{h,H}}_{L^2(D)}^2$. Both region boundaries have zero mass under the limiting measure in~\eqref{eq:common-local}, so almost surely
\begin{equation*}
 \frac{\int_{D_L}|\widetilde L_{h,H}(x)|^2\,dx}
      {\int_D|\widetilde L_{h,H}(x)|^2\,dx}
 \longrightarrow
 \frac{\int_{D_L}\mu(x)\,dx}{\int_D\mu(x)\,dx}.
\end{equation*}
The predicted fractions are $1/2$ for $\mu_0$ and $1/2+1/(2\pi)$ for $\mu_1$. At $r=8$ and $h=1/128$, the measured fractions are $50.61\%$ and $66.54\%$. Their increase of $15.93$ percentage points is close to the predicted $15.92$ percentage points. The whole-region ratios $S_{h,r}/C_\mu$ at the same resolution are $1.00692$ and $1.01390$, with shared fixed-$r$ limit approximately $0.99642$. The two examples illustrate the effect of changing the spatial distribution of the source strength while keeping its integral fixed.

\subsection{Exterior solution and computational cost}

These checks and timings use the constant source strength $\mu=1$. The weight quadrature and FFT application were checked by tighter quadrature, direct radial integration on selected offsets, and direct convolution on a small grid. Appendix~\ref{app:quadrature} gives the quadrature settings and the measured changes in the weights and the common-grid statistic.

For a real source, the outgoing and incoming solutions are complex conjugates and give the same squared norms for differences between successive approximations. To check the outgoing solution, we take $k=4$ and $\alpha=15/16$ and use the cell integrals of the same three-dimensional realization on the fixed source mesh $h=1/8$. We average the field over four fixed cubes of side $h$ outside $D$. Direct quadrature of the outgoing Helmholtz term and the radial Yukawa integral over the source and receiver cells gives the reference values. Refining this quadrature changed the reference vector by $7.8\times10^{-10}$ in relative Euclidean norm. In Table~\ref{tab:receivers}, each halving of $H$ reduces the errors in both components by a factor of about four.

\begin{table}[!htb]
\centering
\caption{Relative Euclidean errors (\%) in the vector of four receiver averages and in its real and imaginary parts.}
\label{tab:receivers}
\small
\begin{tabular}{@{}rrrr@{}}
\toprule
$H$ & Complex field & Real part & Imaginary part \\
\midrule
$1/8$  & 2.060   & 2.050   & 2.108   \\
$1/16$ & 0.5103  & 0.5078  & 0.5224  \\
$1/32$ & 0.1273  & 0.1267  & 0.1303  \\
$1/64$ & 0.03180 & 0.03165 & 0.03256 \\
\bottomrule
\end{tabular}
\end{table}

The computations used Python 3.12.14 and NumPy 2.3.5~\cite{Harris} on a Windows 11 system with an AMD Ryzen AI 7 H 350 processor. We measured the time required to compute the difference between successive numerical solutions on a common output mesh. The wall-clock time includes forming the difference of the discretized sources, generating and embedding the weights, and applying the FFT convolution. The quadrature checks and the evaluation of the error norms were carried out separately. Each case was run three times with weights generated afresh. The values of $S_{h,r}$ agreed across the three runs to the reported precision.

\begin{table}[!htb]
\centering
\caption{Computational time in seconds for the difference between successive numerical solutions on a common output mesh. Each row reports the run with the median total time among three runs, with weights recomputed for each run. The padded grid has $M^d$ points. Setup includes kernel embedding, source assignment, and forward transforms. Apply includes products, inverse transforms, and restriction to $Q$. The total also includes weight generation and the construction of the difference of the discretized sources.}
\label{tab:timing}
\small
\begin{tabular}{@{}rrrrrrrr@{}}
\toprule
$d$ & $h$ & $r$ & $M$ & Weights & Setup & Apply & Total \\
\midrule
2 & $1/128$ & 4 & 2560 & 1.011 & 0.312 & 0.218 & 1.564 \\
2 & $1/128$ & 8 & 5120 & 3.596 & 1.209 & 1.006 & 5.905 \\
3 & $1/16$ & 2 & 160 & 0.129 & 0.203 & 0.159 & 0.500 \\
3 & $1/16$ & 4 & 320 & 0.916 & 1.547 & 1.199 & 3.725 \\
3 & $1/16$ & 8 & 640 & 7.363 & 16.140 & 12.602 & 36.608 \\
3 & $1/32$ & 2 & 320 & 0.950 & 1.590 & 1.250 & 3.855 \\
3 & $1/32$ & 4 & 640 & 7.412 & 16.549 & 12.580 & 37.069 \\
\bottomrule
\end{tabular}
\end{table}

Table~\ref{tab:timing} gives the computational cost for several output mesh sizes at a fixed source mesh size. In three dimensions at $h=1/16$, taking $r=8$ gives a relative difference of about $0.46\%$ between $I_\mu\kappa_r$ and $C_\mu$. The median total time is $36.6$ seconds, including $7.36$ seconds for the weights. At $h=1/32$ and $r=4$, the total times ranged from $36.8$ to $47.8$ seconds. In two dimensions the weights account for about $61\%$ of the median execution time at the finest output resolution.

\section{Conclusion}

The piecewise constant approximation of the white noise gives the conditional mean of the outgoing solution. Its mean square error and its almost-sure squared $L^2$ error have the same leading constant along the dyadic meshes. This constant is optimal among approximations based on the cell integrals, in the sense of Theorem~\ref{thm:path}. The measures defined by the rescaled squared error of the conditional mean converge weakly almost surely to a measure with density proportional to the source strength on $D$. On aligned boxes, the FFT method attains the optimal constant as the ratio of the output mesh size to the source mesh size tends to zero. The numerical examples illustrate the error estimates for the two output mesh choices and the dependence of the spatial error distribution on the source strength.

\appendix
\section{The integrability threshold and the classical endpoint}

Let $0<\alpha\le d/4$ and retain the source assumptions of Section~\ref{sec:model}. The heat representation~\eqref{eq:heat-split} gives the bound $C(1+\abs{x}^{2\alpha-d})$ near the origin, using $\omega(s)\le Cs^{\alpha-1}$. In particular, $G$ is locally integrable. Its high-frequency part is square integrable away from the origin by integration by parts in the Fourier integral, since sufficiently high derivatives of the high-frequency multiplier are integrable. If $G$ were square integrable near the origin, subtracting the smooth low-frequency part would make the high-frequency part globally square integrable. Plancherel's theorem rules this out because the multiplier at high frequencies is comparable to $\abs{\xi}^{-2\alpha}$.

The local $L^1$ bound for $G$ and boundedness of $a$ give a bounded operator $T_\mu$ from $L^2(D)$ to $L^2(Q)$. For every $y\in D$, the set $Q-y$ contains a fixed ball about the origin, since $\overline D\subset Q$. It follows from $I_\mu>0$ that
\begin{equation*}
  \int_D\mu(y)\int_Q\abs{G(x-y)}^2\,dx\,dy=\infty.
\end{equation*}
Thus $T_\mu$ is not Hilbert--Schmidt. Take a real orthonormal basis of the complex output space consisting of smooth compactly supported test functions, and let $\Sigma_n$ be the sum of the first $n$ squared real coefficients of the distributional solution. The expectations $\E\Sigma_n$ tend to infinity because $T_\mu$ is not Hilbert--Schmidt. The covariance operator of these coefficients has norm at most $\norm{T_\mu}_{\op}^2$. Hence~\eqref{eq:variance} gives $\operatorname{Var}(\Sigma_n)\le 2\norm{T_\mu}_{\op}^2\E\Sigma_n$. Chebyshev's inequality yields divergence in probability, which by monotonicity gives divergence almost surely. An $L^2(Q)$ representative would have finite coefficient sum. The same obstruction applies to the residual relative to the finite-rank approximation $T_\mu P_h W\in L^2(Q)$.

For $\alpha=1$ and $d=3$, the Riesz argument of Section~\ref{sec:sharp} still applies, with $p=1$. For $\alpha=1$ and $d=2$, the minimum mean square error on aligned boxes is
\begin{equation}\label{eq:endpoint}
 R_h\sim \frac{I_\mu}{24\pi}h^2\log(1/h).
\end{equation}
As in the proof of Theorem~\ref{thm:risk}, we replace the source coefficient $a$ by $a_h=P_h a$. Local square integrability of $G$ gives an $O(h)$ change in the Hilbert--Schmidt norm of the residual operator. Split the Fourier multiplier of the Green function into a low-frequency part and a high-frequency part supported at $\abs{\xi}>R_0>2k$. The low-frequency part, which is smooth in physical space, contributes $O(h)$ in Hilbert--Schmidt norm. Integrating the high-frequency contribution first over all output space and restricting back to $Q$ changes the squared norm by $O(h^2)$. Indeed, the high-frequency kernel has a square-integrable gradient at a fixed positive distance from the diagonal, and cell Poincar\'e bounds the contribution outside $Q$ by $Ch^2\|a_h\|_{L^2(D)}^2$.

For a source cell $E$, the variance identity gives the whole-space squared norm as $a_E^2|E|/(2\pi)^2$ times the integral of the squared high-frequency multiplier against the cell factor. For $R_0<\abs{\xi}<h^{-1}$ this factor has the expansion
\begin{equation*}
 1-\prod_\ell\sinc^2(h\xi_\ell/2)=h^2\abs{\xi}^2/12+O(h^4\abs{\xi}^4).
\end{equation*}
The squared multiplier equals $\abs{\xi}^{-4}+O(\abs{\xi}^{-6})$ beyond the cutoff transition. The leading radial integral is $2\pi h^2\int_{R_0}^{1/h}\rho^{-1}\,d\rho/12$. The cutoff transition, remainders, and frequency tail contribute $O(h^2)$ to the integral. Since $\sum_E a_E^2\abs{E}\to I_\mu$, the whole-space contribution has the asymptotic~\eqref{eq:endpoint}. The $O(h)$ norm perturbations are negligible at the scale $h\sqrt{\log(1/h)}$. These orders agree with the classical upper estimates in~\cite{CZZ} for constant source strength.

\section{Quadrature for the numerical experiments}\label{app:quadrature}

\subsection{Fixed-ratio limits}

We evaluate the constants $B_r$ in the two fixed-ratio limits by quadrature on a single Fourier cell. For $v\in C_0$ and $\nu\in\mathbb Z^d$, set
\begin{equation*}
 w_\nu(v)=\prod_{\ell=1}^d\sinc^2\bigl(\pi(v_\ell+\nu_\ell)\bigr),\qquad
 a_\nu(v)=\abs{v+\nu}^{-2\alpha},\qquad
 \Theta_r(v)=\prod_{\ell=1}^d\frac{\sinc^2(\pi r v_\ell)}{\sinc^2(\pi v_\ell)} .
\end{equation*}
By the Fourier projection formula for principal shift-invariant spaces~\cite[Theorem~2.9]{BDR}, the output projection at the grouped frequencies $2\pi r(v+\nu)$ is the orthogonal projection onto the vector of Fourier coefficients of a translated output cell indicator. Its squared component magnitudes after normalization are $w_\nu(v)$, whose sum is one by Parseval's identity on $C_0$. The indicator of $C_0$ is the sum of its output cell indicators. At the grouped frequencies their relative translation phases depend only on $v$, so summing the geometric progressions gives the common factor with squared magnitude $\Theta_r(v)$. The Fourier coefficients of $J$ also have the multiplier factor $(2\pi r)^{-2\alpha}a_\nu(v)$. Thus
\begin{equation}\label{eq:b-alias}
 B_r=(2\pi)^{-4\alpha}r^{-p}\int_{C_0}\Theta_r(v)
 \left(\sum_\nu w_\nu(v)a_\nu(v)^2-\left(\sum_\nu w_\nu(v)a_\nu(v)\right)^2\right)\,dv .
\end{equation}
The formula follows first with a low-frequency cutoff on $J$, which can be removed by Poincar\'e's inequality and convergence in the gradient norm. For an integer cutoff $M_{\rm a}\ge0$, we retain the terms with $\max_\ell\abs{\nu_\ell}\le M_{\rm a}$. The omitted profile at frequencies with at least one coordinate of magnitude exceeding $R=2\pi r(M_{\rm a}+1/2)$ has squared $L^2$ norm at most $(4d/\pi)R^{-4\alpha-1}$. This follows by integrating the corresponding one-dimensional sinc tail and applying Parseval to the other coordinates. Its norm bounds the truncation error in $\sqrt{B_r}$.

For the common-grid ratios in Table~\ref{tab:accuracy}, the bound on the contribution of the omitted modes is at most $4.25\times10^{-4}$. The estimated ratios change by at most $1.31\times10^{-5}$ when the Gauss order and the frequency cutoff are refined.

\subsection{Weight evaluation and numerical checks}

For both trapezoidal rules in logarithmic variables, we bound the quadrature and truncation errors separately. In exact arithmetic, the prescribed absolute error tolerance for the heat contribution to each weight is $10^{-7}H$. For broad Gaussians with $s/H^2\ge16$, the one-dimensional average is evaluated by a Gauss--Legendre rule~\cite[Section~3.5(v)]{DLMF} with at least 12 points on each half of the tent support to avoid cancellation in~\eqref{eq:heat-average}. For $\lambda=H^2/(4s)\le1/64$, the function $e^{-\lambda(j+z)^2}$ is bounded by $e$ on the complex disk $\abs z\le8$, uniformly in the one-dimensional offset $j$. The Taylor tail from degree 23 onward is at most $e\,8^{-23}/(1-1/8)$ on $\abs z\le1$. Multiplication of the degree-22 Taylor polynomial by the linear tent factor stays within the exactness range of the rule. The remaining Lipschitz part is evaluated at the displacement center as in Theorem~\ref{thm:weights}.

On the finest output lattices in two and three dimensions, the maximum changes in the heat contribution to the convolution weights after reducing the tolerance by a factor of $1000$ were $3.98\times 10^{-12}$ and $2.80\times 10^{-11}$. The original tolerances were $4.88\times 10^{-11}$ and $3.91\times 10^{-10}$. Independent radial Yukawa integration on selected offsets gave discrepancies below $3.0\times 10^{-12}$ in both dimensions. For the Lipschitz term, the largest midpoint deviations over the sampled offsets from tensor Gauss integration were $7.72\times 10^{-8}$ and $2.19\times 10^{-4}$. These are within the uniform midpoint bounds $1.504H$ and $0.08489H$. Direct summation on a small mesh agreed with the FFT to within $7.0\times10^{-17}$.

We also recomputed all the weights with a smaller tolerance for the heat quadrature and Gaussian integration of the Lipschitz term at every offset. The resulting change in $S_{h,r}/C_\mu$ was $5.55\times10^{-6}$ for the two-dimensional case $(h,r)=(1/16,4)$ and $2.58\times10^{-4}$ for the three-dimensional case $(h,r)=(1/8,2)$.

\section*{Declaration of AI Use}
The authors used OpenAI GPT 6.2 to assist with the writing and presentation of the
manuscript and the development of the code. The authors reviewed the entire
manuscript and all the code and take full responsibility for their content.

\section*{Acknowledgements}
X. Xu is partially supported by the National Key Research and Development Program of China (No. 2024YFA1012303), Fundamental and Interdisciplinary Disciplines Breakthrough Plan of the Ministry of Education of China (JYB2025XDXM805), National Natural Science Foundation of China (No. 12525112), and Zhejiang Key Laboratory of Airspace Awareness and Autonomous Unmanned Systems. Y. Zhao is supported by National Natural Science Foundation of China (No. 12371421).

\section*{Conflict of Interest Statement}
The authors declare that they have no conflict of interest.

\section*{Data Availability Statement}
The numerical experiments use synthetic data generated from the models and parameters described in this article. The accompanying reproducibility archive contains the scripts, CSV results, and instructions for regenerating the figures and tables.

\end{document}